\documentclass[12pt]{article}
\usepackage{latexsym, amsfonts, amsmath}

\newtheorem{thm}{Theorem}[section]
\newtheorem{cor}[thm]{Corollary}
\newtheorem{lem}[thm]{Lemma}

\def\infint{\int_{-\infty}^\infty}

\def\convd{\stackrel{\cal D}{\rightarrow}}

\def\ex{{\rm E\,}}

\def\var{\mathop{\rm Var}\nolimits}

\begin{document}

\bibliographystyle{plain}

\title {Asymptotic Normality of Kernel Type\\
Deconvolution Estimators}

\author {A.J. van Es \quad {\normalsize and}\quad H.-W. Uh\\[.3cm]
{\normalsize Korteweg-de Vries Institute for Mathematics,
 University of Amsterdam}\\
{\normalsize Plantage Muidergracht 24,
 1018 TV Amsterdam}\\
{\normalsize The Netherlands}}
\date{}
\maketitle

\begin{abstract}
We derive asymptotic normality of kernel type deconvolution estimators
of the density, the distribution function at a fixed point, and of the
probability of an interval.  We consider the so called
super smooth case where the characteristic function of the known
distribution decreases exponentially.

It turns out that the limit behavior of the pointwise estimators of
the density and distribution function is relatively straightforward
while the asymptotics of the estimator of the probability of an
interval depends in a complicated way on the sequence of bandwidths.
\\[.5cm] 

\end{abstract}

\section{Introduction and results}

Let $X_1,\ldots, X_n$ be i.i.d. observations, where $ X_i=Y_i+Z_i $
and $Y_i$ and $Z_i$ are independent.  Assume that the unobservable
$Y_i$ have distribution function $F$ and density $f$. Also assume that
the random variables $Z_i$ have a known density $k$. Note that the
density $g$ of $X_i$ is equal to the convolution of $f$ and $k$.  The
deconvolution problem is the problem of estimating $f$ or $F$ from the
observations $X_i$.

Kernel type estimators for the density $f$ and its distribution
function $F$ have been studied by many authors. We mention Carroll and
Hall (1988), Liu and Taylor (1989), Stefanski (1990), Stefanski and
Carroll (1990), Zhang (1990), Fan (1991a,b, 1992), Fan and Liu (1997)
and Van Es and Kok (1998).

Let $w$ denote a {\em kernel function} and $h>0$ a {\em bandwidth}.
The estimator $f_{nh}(x)$ of the density $f$ at the point $x$ is defined as
\begin{equation}\label{fnh}
f_{nh}(x)=\frac{1}{2\pi} \int_{-\infty}^\infty
e^{-itx} \frac{\phi_w(ht)\phi_{emp}(t)}{ \phi_k(t)}\,dt.
\end{equation}
Here $\phi_{emp}$ denotes the empirical
characteristic function of the sample, i.e.
$$%\begin{equation}
\phi_{emp}(t) = {1\over n}\sum_{j=1}^n e^{itX_j},
$$%\end{equation}
and $\phi_w$ and $\phi_k$ denote the characteristic functions of
$w$ and $k$ respectively.
A straightforward estimator   of
$P(a<Y_j\leq b)
=F(b)-F(a)$  is, for $-\infty<a<b<\infty$, given by
\begin{equation}\label{Fnh}
F_{nh}(a,b)=\int_a^b f_{nh}(x).
\end{equation}
Computing the expectation of $f_{nh}(x)$ we get
\begin{eqnarray*}
\ex f_{nh}(x)&=& \ex \frac{1}{2\pi}\int_{-\infty}^\infty
e^{-itx}\frac{\phi_w(ht)\phi_{emp}(t)}{\phi_k(t)}dt\\
&=&
\frac{1}{2\pi}\int_{-\infty}^\infty
e^{-itx}\frac{\phi_w(ht)\phi_g(t)}{\phi_k(t)}dt\\
&=&
\frac{1}{2\pi}\int_{-\infty}^\infty
e^{-itx}\phi_w(ht)\phi_f(t)\,dt\\
&=&
\int_{-\infty}^\infty \frac{1}{h} w\Big(\frac{x-u}{h}\Big)f(u)\,du,
\end{eqnarray*}
which shows that this expectation is equal to the expectation of an
ordinary kernel estimator based on observations from $f$.

Most of the papers on kernel type deconvolution deal with mean squared
error properties and optimal rates of convergence. The rate of decay
to zero, at minus infinity and plus infinity, of the modulus of the
characteristic function $\phi_k$ is crucial to the asymptotics. Two
cases have been distinguished, the {\em smooth case}, where $|\phi_k|$
decays algebraically to zero, and the {\em super smooth case}, where
it decreases exponentially.

Asymptotic normality of the density estimator $f_{nh}(x)$ has been
derived by Zhang (1990), Fan (1991b) and Fan and Liu (1997).  Zhang
considers asymptotic normality of the estimator of the distribution
function as well.  Asymptotic normality of the density estimator based
on a stationary sequence of observations has been established by Masry
(1993).  In the smooth case the asymptotics are essentially the same
as those of higher order derivatives of an ordinary kernel estimator
of $g(x)$. See Van Es and Kok (1998) for some relatively simple
deconvolution problems, such as gamma and Laplace deconvolution, where
this is obvious. Here the limit variance depends on the (unknown)
value of $f(x)$.  In the super smooth case the asymptotics are much
more complicated. Asymptotic normality has been established for
studentized estimators, i.e. the difference between the estimators and
their expectation is divided by an estimate of the standard
deviation. The asymptotic behavior of the variance however is not
clear.  This random standardization is also motivated by the need of
confidence intervals.

In Section \ref{sect1.1} of this paper for simplicity we consider normal
deconvolution where $k$ is the standard normal density, so
\begin{equation}\label{normphi}%$$
\phi_k(t)=e^{-\frac{1}{2}t^2}.
\end{equation}%$$
The results are generalized in Section \ref{sect1.2} where the general
super smooth case is considered.  It turns out that we have to put a
restriction on the rate of decay at plus and minus infinity of
$\phi_k$.  Essentially only rates faster than $e^{-|t|}$ are allowed,
thus excluding for instance the Cauchy distribution.

Our results below indicate that the asymptotic variance in the cases
considered in this paper is {\em independent of $f$ and $x$}, so
random standardization to obtain confidence intervals is not
necessary. It turns out that while the asymptotic behavior of the
density estimator (\ref{fnh}) and a suitable estimator of $F(b)$,
established in Theorem \ref{supsmodens} and Theorem \ref{dstrb} below,
are relatively straightforward the asymptotic behavior of the
estimator (\ref{Fnh}) of $F(b)-F(a)$, in Theorem \ref{supsmodist}, is
rather unusual and complicated as it depends heavily on the sequence
of bandwidths.

\bigskip

Throughout the paper we impose the following condition on the kernel
function $w$.

\bigskip

\noindent{\bf Condition W}

Let $\phi_w$ be real valued, symmetric and have support $[-1,1]$.  Let
$\phi_w(0)=1$, and let
\begin{equation}\label{condfiw}
\phi_w(1-t)=At^\alpha + o(t^\alpha), \qquad \mbox{as} \qquad t\downarrow 0
\end{equation}
for some constants $A$ and $\alpha\geq 0$.

\medskip

The remainder of the paper is organized as follows.  In Section
\ref{sect1.1} we give the results for normal deconvolution, and in
Section \ref{sect1.2} for general super smooth deconvolution.  In
Section \ref{theoremproofs} we give the main steps of the proofs of
the three theorems while in Section \ref{lemmaproofs} we prove the
lemmas containing the details.

\subsection{Normal deconvolution}\label{sect1.1}

To simplify the presentation we consider normal deconvolution first,
i.e. we assume (\ref{normphi}).  Our first theorem establishes
asymptotic normality of the density estimator (\ref{fnh}), for a
sequence of bandwidths $h_n$, i.e. we consider $f_{nh_n}$.  The
dependence of the bandwidth $h_n$ on $n$ is suppressed in most of the
sequel.

\begin{thm}\label{supsmodens}

Assume Condition W and   $\ex X^2<\infty$.
Then, as $n\to\infty$ and $h\to 0$,
$$
\frac{\sqrt{n}}{h^{1+2\alpha}e^{\frac{1}{2h^2}}}\,
(f_{nh}(x) - \ex f_{nh}(x))\convd 
N(0,\frac{A^2}{2\pi^2}\Gamma(\alpha+1)^2),
$$
where $\Gamma(t)=\int_0^{\infty}v^{t-1} e^{-v}dv$.

\end{thm}

To establish asymptotic normality of the estimator (\ref{Fnh}) we need
an extra, rather complicated, condition on the sequence of bandwidths
$h$.  Consider $u=1/(2h)$ instead of $h$.  Note that as $h \rightarrow
0$ we have $u\rightarrow \infty$.  Let $\mathcal{S}$ denote the set of
positive points $u$ where $\sin((b-a)u)$ vanishes, i.e.
$$
\mathcal{S}=\{ \frac{\pi k}{b-a}\,: \,k=0,1,2,\ldots\}.
$$
Denote by $u^-$ the largest element of $\mathcal{S}$ below $u$, and by
$u^+$ the smallest element above $u$.  If $u$ belongs to $\mathcal{S}$
then $u, u^-$ and $u^+$ coincide. The following condition
specifies how close $u$ is to $\mathcal{S}$. Surprisingly we get
different asymptotic behavior of $F_n(a,b)$ in these cases.

\bigskip

\noindent{\bf Condition A}

As $n \rightarrow \infty$ we have for $u=u_n$
\begin{description}
\item[A1:]
$u(u-u^-)\rightarrow \infty $ and $u(u^+-u)\rightarrow
\infty,$
\item[A2:] 
$\min(u(u-u^-), u(u^+-u))\rightarrow
0,$
\item[A3:]
for some constant $\gamma>0$ either
$u(u-u^-)\rightarrow \gamma$ or $u(u^+-u)\rightarrow
\gamma$.  
\end{description}

\begin{thm}\label{supsmodist}
Assume Condition W and  $\ex X^2<\infty$.
Then, as $n\to\infty$ and $h\to 0$,\\
(a) under Condition A1,
$$
\frac{\sqrt{n}}{h^{2+2\alpha}e^{\frac{1}{2h^2}}\sin(\frac{b-a}{2h})}\,
(F_{nh}(a,b) - \ex F_{nh}(a,b))\convd 
N\Big(0,\, \frac{2A^2}{\pi^2}\Gamma(\alpha+1)^2\Big),
$$
(b) under Condition A2,
$$%\begin{equation}
\frac{\sqrt{n}}{h^{3+2\alpha}e^{\frac{1}{2h^2}}}\,
(F_{nh}(a,b) - \ex F_{nh}(a,b))\convd N\Big(0,\,\frac{A^2}{2\pi^2}
\Gamma(\alpha+2)^2(b-a)^2\Big),
$$%\end{equation}
(c) under Condition A3,
$$%\begin{equation}
\frac{\sqrt{n}}{h^{3+2\alpha}e^{\frac{1}{2h^2}}}\,
(F_{nh}(a,b) - \ex F_{nh}(a,b))\convd N\Big(0,\,\frac{A^2}{2\pi^2}
(4\gamma\Gamma(\alpha+1)+\Gamma(\alpha+2))^2(b-a)^2\Big).
$$%\end{equation}
\end{thm}

The next corollary gives an order bound for the four different cases
of this theorem.

\begin{cor}
Under the conditions of Theorem \ref{supsmodist} we have
$$
F_{nh}(a,b) - \ex F_{nh}(a,b)=O_P\Big(
\frac{1}{\sqrt{n}}\,h^{2+2\alpha}e^{\frac{1}{2h^2}}\Big).
$$
\end{cor}

Finally we consider the asymptotics of an estimator of $F(b)$.
The next theorem establishes asymptotic normality of $F_{nh}(a,b)$
with $a$ tending to minus infinity at a suitable rate.

\begin{thm}\label{dstrb}
Assume Condition W and $E X^2 < \infty$.  Let $n\to \infty, h\to 0$
and $a\to-\infty$ such that $ah\rightarrow -\infty$ and
$a=o(e^{(1/2h^2)(1-\delta)})$ for some $0<\delta<1$.  Then, for $b$
fixed, we have
$$
\frac{\sqrt{n}}{h^{2+2\alpha}e^{\frac{1}{2h^2}}}\,
(F_{nh}(a,b) - \ex F_{nh}(a,b))\convd N\Big(0,\,\frac{A^2}{2\pi^2} 
\Gamma(\alpha+1)^2\Big).
$$
\end{thm}

\subsection{General super smooth deconvolution}\label{sect1.2}

In this section we generalize the results for normal deconvolution to
the super smooth case, i.e. we consider densities $k$ that satisfies
the following condition.

\bigskip

\noindent{\bf Condition K}

Assume that $\phi_k$ satisfies
as follows.
$$
\phi_k(t)\sim C |t|^{\lambda_0} e^{-|t|^\lambda/\mu},
$$
as $|t|\rightarrow \infty$ for some $\lambda>1, \mu>0, \lambda_0$, and
some real constant $C$.  Furthermore assume $\phi_k(t) \neq 0$ for all
$t$.

\medskip

Note that this condition excludes the Cauchy distribution and all
other distributions for which the tail of the characteristic function
decreases slower than $e^{-|t|}$.

\begin{thm}\label{supsmodens*}
Assume Condition W, Condition K and   $\ex X^2<\infty$.
Then, as $n\to\infty$ and $h\to 0$,
$$
\frac{\sqrt{n}}{h^{\lambda(1+\alpha)+\lambda_0-1}
e^{\frac{1}{\mu h^\lambda}}}\,(f_{nh}(x)-\ex f_{nh}(x))\convd 
N(0,\frac{A^2}{2\pi^2}(\mu/ \lambda)^{2+2\alpha}(\Gamma(\alpha+1))^2).
$$
\end{thm}

\bigskip

\noindent We need a slightly different condition on the sequence of
the bandwidths $h$ to generalize Theorem \ref{supsmodist}.%th1.2

\bigskip

\noindent{\bf Condition A*}

As $n \rightarrow \infty$ we have for $u=u_n$
\begin{description}
\item[A1*:]
$u^{\lambda-1}(u-u^-)\rightarrow \infty $ and $u^{\lambda-1}(u^+-u)\rightarrow
\infty,$
\item[A2*:] 
$\min(u^{\lambda-1}(u-u^-), u^{\lambda-1}(u^+-u))\rightarrow
0,$
\item[A3*:]
for some constant $\gamma>0$ either
$u^{\lambda-1}(u-u^-)\rightarrow \gamma$ or $u^{\lambda-1}(u^+-u)\rightarrow
\gamma$.
\end{description}

\medskip

\begin{thm}\label{supsmodist*}
Assume Condition W and  $\ex X^2<\infty$.
Then, as $n\to\infty$ and $h\to 0$,\\
(a) under Condition A1*,
$$
\frac{\sqrt{n}}
{h^{(1+\alpha)\lambda+\lambda_0}e^{\frac{1}{\mu h^\lambda}}\sin(\frac{b-a}{2h})}\,
(F_{nh}(a,b) - \ex F_{nh}(a,b))\convd 
N\Big(0,\, \frac{2A^2}{\pi^2}(\mu/ \lambda)^{2+2\alpha}
\Gamma(\alpha+1)^2\Big),
$$\\
(b) under Condition A2*,
$$
\frac{\sqrt{n}}
{ h^{(2+\alpha)\lambda+\lambda_0-1}e^{\frac{1}{\mu h^\lambda}} }\,
(F_{nh}(a,b) - \ex F_{nh}(a,b))\convd N\Big(0,\,\frac{A^2}{2\pi^2}
(\mu /\lambda)^{4+2\alpha}\Gamma(\alpha+2)^2(b-a)^2\Big),
$$
(c) under Condition A3* ,
$$
\frac{\sqrt{n}}{h^{(2+\alpha)\lambda+\lambda_0-1}
e^{\frac{1}{\mu h^\lambda}}}\,
(F_{nh}(a,b) - \ex F_{nh}(a,b))\convd 
N\Big(0,\,\frac{A^2}{2\pi^2}
\Big(2^\lambda \gamma \Gamma(\alpha+1)+ (\mu/\lambda)\Gamma(\alpha+2)\Big)^2
(\mu/\lambda)^{2+2\alpha}(b-a)^2\Big).
$$
\end{thm}

\noindent The next corollary gives an order bound for the four
different cases of this theorem.
\begin{cor}
Under the conditions of Theorem \ref{supsmodist*} we have
$$
F_{nh}(a,b) - \ex F_{nh}(a,b)=O_P\Big(
\frac{1}{\sqrt{n}}\,h^{(1+\alpha)\lambda+\lambda_0}e^{\frac{1}{\mu h^\lambda}}\Big).
$$
\end{cor}

\noindent Theorem \ref{dstrb} is generalized as follows.
\begin{thm}\label{dstrb*}
Assume Condition W and Condition K.  Suppose $E X^2 < \infty$.  Let
$n\to \infty, h\to 0$ and $a\to-\infty$ such that
$ah^{\lambda-1}\rightarrow -\infty$ and $a=o(e^{(1/\mu
h^\lambda)(1-\delta)})$ for some $0<\delta<1$.  Then, for $b$ fixed,
we have
$$
\frac{\sqrt{n}}{h^{(1+\alpha)\lambda+\lambda_0}
e^{\frac{1}{\mu h^\lambda}} }\,
(F_{nh}(a,b) - \ex F_{nh}(a,b))\convd 
N\Big(0,\,\frac{A^2}{2\pi^2} (\mu/\lambda)^{2+2\alpha}
\Gamma(\alpha+1)^2\Big).
$$
\end{thm}

\bigskip
Clearly the arguments in the proofs fail for $\lambda \leq 1$.  To see
if this is essential we consider the deconvolution density estimator
for Cauchy deconvolution, where $\phi_k(t)=e^{-|t|}$, with a special
choice of $\phi_w$.  For $\phi_w(t)=I_{[-1,1]}(t)$ we have $\alpha=0,
A=1$ and
\begin{eqnarray*}
f_{nh}(x)&=& \frac{1}{2\pi n}\sum_{j=1}^n \int_{-1/h}^{1/h}
\cos(t(X_j-x))e^{|t|}dt\\
&=&
\frac{1}{\pi n} \sum_{j=1}^n
\frac{1}{1+(X_j-x)^2}
\Big\{-1+e^{1/h}\Big(\cos\Big(\frac{X_j-x}{h}\Big)+
(X_j-x)\sin\Big(\frac{X_j-x}{h}\Big)\Big)\Big\}.
\end{eqnarray*}
Asymptotic normality of $f_{nh}$ is given by the following theorem.
The proof is given in Section \ref{s24}.

\begin{thm}\label{cauchy}
As $n\to \infty$ and $h\to 0$ we have
%\begin{gather*}
$$
\sqrt{n} e^{-1/h}(f_{nh}(x)-\ex f_{nh}(x))\convd 
N(1, \sigma^2), $$%\\\text{with}\qquad 
with
$$
\sigma^2=\frac{1}{2\pi^2} \ex \frac{1}{1+(X_j-x)^2}
= \frac{1}{2\pi^2} \int_{-\infty}^\infty \frac{1}{1+(u-x)^2}g(u)\,du.
$$%\end{gather*}
\end{thm}
It follows that the rate of convergence is equal to the one in Theorem
\ref{supsmodens*} for $\alpha=0, \lambda=1$ and $\lambda_0=0$.
However the limit variance now depends on $x$, and on $f$ through
$g=f*k$.  This shows that the condition $\lambda>1$ is essential.

\section{Proofs of the theorems}\label{theoremproofs}
\subsection{Proofs of Section \ref{sect1.1}}

\subsubsection{Proof of Theorem \ref{supsmodens}}

The theorem is immediate from Lemma \ref{asnor} in Section
\ref{lemmaproofs}, and the following lemma that represents $f_{nh}(x)$
asymptotically as a normalized mean.
\begin{lem}\label{decomp1}
$$
\frac{\sqrt{n}}{h^{1+2\alpha}e^{\frac{1}{2h^2}}}\,
(f_{nh}(x)-\ex f_{nh}(x))
=
\frac{A}{2 \pi}(\Gamma(\alpha+1)+o(1))U_{nh}(x) +O_P(h),
$$
where
$$
U_{nh}(x) ={1\over \sqrt{n}}\sum_{j=1}^n
\Big(\cos\Big(\frac{X_j-x}{h}\Big)-\ex \cos\Big(\frac{X_j-x}{h}\Big)\Big).
$$
\end{lem}

\subsubsection{Proof of Theorem \ref{supsmodist}}%th1.2

We can represent $F_n(a,b)$ asymptotically as a constant $\tau_n$
times a normalized mean.

\begin{lem}\label{decomp2}
\begin{eqnarray*}
\lefteqn{
\frac{\sqrt{n}}{h^{2+2\alpha}e^{\frac{1}{2h^2}}}\,
(F_{nh}(a,b)-\ex F_{nh}(a,b))=\tau_n S_{nh}(a,b)}\\
&+&
O_P\Big(\frac{1}{h^{3+2\alpha}}\,e^{\frac{1}{2h^2}
(\epsilon^2-1)}\Big)
+
O_P\Big(h\Big|\sin\Big(\frac{b-a}{2h}\Big)\Big| \Big) +
O_P(h^2),
\end{eqnarray*}
where
\begin{equation}\label{taundef}
\tau_n=\frac{A}{\pi}
\Big(
2(\Gamma(\alpha+1)+o(1))\sin\Big(\frac{b-a}{2h}\Big)
+(b-a)(\Gamma(\alpha+2)+o(1))\cos\Big(\frac{b-a}{2h}\Big)h
\Big)
\end{equation}
and
\begin{equation}\label{Sndef}
S_{nh}(a,b)=\frac{1}{\sqrt{n}}\sum_{j=1}^n \Big(\cos\Big(\frac{X_j-(a+b)/2}{h}
\Big)
-\ex \cos\Big(\frac{X_j-(a+b)/2}{h}\Big)\Big).
\end{equation}
\end{lem}

\bigskip

By Lemma \ref{asnor} with $x$ equal to $(a+b)/2$ we have asymptotic
normality of $S_{nh}(a,b)$, i.e.
$$
S_{nh}(a,b) \convd N(0, \frac{1}{2}).
$$

The proof of the theorem is finished once we establish the asymptotic
behavior of $\tau_n$ under conditions A1, A2 and A3.  Under A1 we
have three possible situations: (i) $u$ remains at a fixed distance of
$\mathcal{S}$, (ii) $u-u^- \rightarrow 0$ and (iii) $u^+-u \rightarrow
0$.  In case (i) we have
$$
\Big|\sin\Big(\frac{b-a}{2h}\Big)\Big|=|\sin((b-a)u)|\gg  \frac{1}{u}=2h.
$$
This also holds in situation (ii), since
\begin{eqnarray*}
\lefteqn{
\Big|\sin\Big(\frac{b-a}{2h}\Big)\Big|=|\sin((b-a)u)- \sin((b-a)u^-)|}\\
&=&
\Big|\frac{\sin((b-a)u)- \sin((b-a)u^-)}{ (b-a)(u-u^-) }\Big| (b-a)(u-u^-)\\
&\sim& 
1\cdot (b-a)(u-u^-) \gg {1\over u}=2h.
\end{eqnarray*}
A similar argument holds in situation (iii).
So Condition A1 implies
$$
\Big|\sin\Big(\frac{b-a}{2h}\Big)\Big| \gg h.
$$
Similarly Condition A2 implies
$$
\Big|\sin\Big(\frac{b-a}{2h}\Big)\Big| =o(h).
$$
To deal with Condition A3 we distinguish two sequences of $\{u_n\}$,
the subsequence $\{u_{n_i}\}$ for which the $k_{n_i}$ corresponding to
the closest points in $\mathcal{S}$ are even, and the subsequence
$\{u_{n_j}\}$ for which the $k$'s are odd.  Along the first
subsequence we have $\cos ((b-a)/2h)\sim 1$, while along the second
subsequence $\cos ((b-a)/2h)\sim -1$.  Moreover we have, along the
first subsequence
$$
\sin\Big(\frac{b-a}{2h}\Big) \sim 2(b-a)\gamma h,
$$
and along the second subsequence
$$
\sin\Big(\frac{b-a}{2h}\Big) \sim -2(b-a)\gamma h.
$$
All this implies
\begin{eqnarray*}
|\tau_n|
&=&{A\over \pi}
\Big(
2(\Gamma(\alpha+1)+o(1))\Big|\sin\Big(\frac{b-a}{2h}\Big)\Big|
+(b-a)(\Gamma(\alpha+2)+o(1))\Big|\cos\Big(\frac{b-a}{2h}\Big)\Big|h
\Big)\\
&\sim&
\begin{cases}
\frac{2A}{\pi}\,\Gamma(\alpha+1)\Big|\sin\Big(\frac{b-a}{2h}\Big)\Big|
&\qquad \text{under A1},\\
\frac{A}{\pi}\,\Gamma(\alpha+2)(b-a)h
&\qquad \text{under A2},\\
\frac{A}{\pi}\,\Big(4\gamma\Gamma(\alpha+1)+\Gamma(\alpha+2)\Big)(b-a)h
&\qquad \text{under A3 }.
\end{cases}
\end{eqnarray*}
By computing the asymptotic variances the proof is then completed.

\subsubsection{Proof of Theorem \ref{dstrb}}

Choose $0<\epsilon< \sqrt{\delta}$.
>From the proof of Lemma \ref{decomp2} we have
$$
F_{nh}(a,b)=
\frac{1}{2\pi n}
\sum_{j=1}^n \int_{-1}^1 {1 \over {is}}
(e^{is\frac{X_j-a}{h}}-e^{is\frac{X_j-b}{h}}) \phi_w(s)
e^{\frac{1}{2h^2}s^2}ds.
$$
We can decompose $F_{nh}(a,b)$ as follows
\begin{eqnarray}
\lefteqn{F_{nh}(a,b)=W_{nh}(b)+ R_{nh}(a)
\nonumber}\\
&&
+\frac{1}{2\pi n}
\sum_{j=1}^n \int_{-\epsilon}^\epsilon \frac{1}{ is}
(e^{is\frac{X_j-a}{h}}-e^{is\frac{X_j-b}{h}}) \phi_w(s)
e^{\frac{1}{2h^2}s^2}ds,\label{term3}
\end{eqnarray}
with
$$
W_{nh}(b)=\frac{1}{\pi n}
\sum_{j=1}^n \int_{\epsilon}^1 \frac{1}{s}
\sin\Big(s\Big(\frac{X_j-b}{h}\Big)\Big) \phi_w(s)
e^{\frac{1}{2h^2}s^2}ds
$$
and
\begin{equation}\label{rnh213}%$$
R_{nh}(a)=\frac{1}{2\pi n}
\sum_{j=1}^n \Big(\int_{-1}^{-\epsilon}+\int_{\epsilon}^1\Big) \frac{1}{is}
e^{is\frac{X_j-a}{h}} \phi_w(s)e^{\frac{1}{2h^2}s^2}ds.
\end{equation}%$$
The next lemma establishes asymptotic normality of $W_{nh}(b)$.
Its proof is omitted because of its similarity to the proofs
of Theorem  \ref{supsmodens} and Theorem \ref{supsmodist}.
\begin{lem}\label{asnor1}
$$
\frac{\sqrt{n}}{h^{2+2\alpha}e^{\frac{1}{2h^2}}}\,
(W_{nh}(b) - \ex W_{nh}(b))\convd N(0,\frac{A^2}{2\pi^2}
\Gamma(\alpha+1)^2).
$$
\end{lem}

For the term $R_{nh}(a)$ we have the following order bound.
\begin{lem}\label{rllem}
\begin{equation}\label{rleq}
R_{nh}(a)-\ex R_{nh}(a)=
o_P\Big(\frac{1}{\sqrt{n}}\,h^{2+2\alpha}e^{\frac{1}{2h^2}}\Big).
\end{equation}
\end{lem}

Since the absolute value of the term (\ref{term3}) is smaller than
$(1/\epsilon \pi h)(b-a)e^{\frac{1}{2h^2}\epsilon^2}$
a combination of the two lemmas proves the theorem.

\subsection{Proofs of Section \ref{sect1.2}}

\subsubsection{Proofs of Theorem \ref{supsmodens*}}

Since the proof is somewhat similar to the proof for normal
deconvolution in the previous section we will state only the
corresponding lemma.
\begin{lem}\label{decomp1*}
\begin{eqnarray*}
\lefteqn{
\frac{\sqrt{n}}{h^{\lambda(1+\alpha)+\lambda_0-1}
e^{\frac{1}{\mu h^\lambda}}}\,(f_{nh}(x)-\ex f_{nh}(x))}\\
&=&
\frac{A}{\pi}(\frac{\mu}{\lambda})^{1+\alpha}
(\Gamma(\alpha+1)+o(1))U_{nh}(x) +O_P(h^{\lambda-1})
+O_P\Big(\epsilon^{1-\lambda_0}h^{-\lambda(1+\alpha)}
e^{ \frac{1}{\mu h^\lambda}(\epsilon^{\lambda}-1)}\Big) ,
\end{eqnarray*}
where
$$
U_{nh}(x) =\frac{1}{\sqrt{n}}\sum_{j=1}^n
\Big(\cos\Big(\frac{X_j-x}{h}\Big)-\ex \cos\Big(
\frac{X_j-x}{h}\Big)\Big).
$$
\end{lem}

\subsubsection{Proofs of Theorem \ref{supsmodist*}}%th1.6

We first write $F_{nh}(a,b)$ asymptotically as $\tau_n$ times a
normalized mean.
\begin{lem}\label{decomp2*}%lem2.6
\begin{eqnarray*}
\lefteqn{\frac{\sqrt{n}}{h^{(1+\alpha)\lambda+\lambda_0}
e^{\frac{1}{\mu h^\lambda}}}\,(F_{nh}(a,b)-\ex F_{nh}(a,b))
=\tau_n S_{nh}(a,b)}\\
&+&
O_P\Big(h^{-(1+\alpha)\lambda-1}\epsilon^{-\lambda_0}\,
e^{\frac{1}{\mu h^\lambda}(\epsilon^\lambda-1)}\Big)
+
O_P\Big(h^{\lambda-1}\Big|\sin\Big({b-a\over 2h}\Big)\Big| \Big) +
O_P\Big(h^{2\lambda-2}\Big),
\end{eqnarray*}
where
$$
\tau_n=\frac{A}{\pi}
\Big(
2\Big(\frac{\mu}{\lambda}\Big)^{1+\alpha}(\Gamma(\alpha+1)+o(1))
\sin\Big({b-a\over 2h}\Big)
+(b-a)\Big(\frac{\mu}{\lambda}\Big)^{2+\alpha}(\Gamma(\alpha+2)+o(1))
\cos\Big(\frac{b-a}{2h}\Big)h^{\lambda-1}
\Big)
$$
and $S_{nh}$ is defined as (\ref{Sndef}).
\end{lem}
Apart from the slightly different condition A* on the sequence of
bandwidths $h$, the proof is similar to the proof for the normal
deconvolution.  Here we have
\begin{equation*}
|\tau_n|\sim
\begin{cases}
\frac{2A}{\pi}\,\Big(\frac{\mu}{\lambda}\Big)^{1+\alpha}\Gamma(\alpha+1)
\Big|\sin\Big(\frac{b-a}{2h}\Big)\Big|
&\qquad \text{under A1*},\\
\frac{A}{\pi}\,\Big(\frac{\mu}{\lambda}\Big)^{2+\alpha}\Gamma(\alpha+2)
(b-a)h^{\lambda-1}
&\qquad \text{under A2*},\\
\frac{A}{\pi}\,\Big(\frac{\mu}{\lambda}\Big)^{1+\alpha}\Big(2^\lambda \gamma\Gamma(\alpha+1)+(\frac{\mu}{\lambda})\Gamma(\alpha+2)\Big)(b-a)h^{\lambda-1}
&\qquad \text{under A3* },
\end{cases}
\end{equation*}
which yields the asymptotic variances of the theorem.

\subsubsection{Proof of Theorem \ref{dstrb*}}

Choose $0<\epsilon<\delta^{1\over \lambda}$.  Let $W_{nh}(b)$ and $
R_{nh}(a)$ be defined by
\begin{equation}
W_{nh}(b)={1\over \pi n}
\sum_{j=1}^n \int_{\epsilon}^1 {1 \over s}
\sin\Big(s\Big({X_j-b\over h}\Big)\Big) \phi_w(s)
\frac{1}{\phi_k(s/h)}ds,
\end{equation}
and
\begin{equation}
R_{nh}(a)={1\over 2\pi n}
\sum_{j=1}^n \Big(\int_{-1}^{-\epsilon}+\int_{\epsilon}^1\Big) {1 \over is}
e^{is{X_j-a\over h}} \phi_w(s)\frac{1}{\phi_k(s/h)}ds.
\end{equation}
Then $F_{nh}(a,b)$ can be written as
\begin{eqnarray}
\lefteqn{F_{nh}(a,b)=W_{nh}(b)+ R_{nh}(a)
\nonumber}\\
&&
+\frac{1}{2\pi n}
\sum_{j=1}^n \int_{-\epsilon}^\epsilon \frac{1}{ is}
(e^{is\frac{X_j-a}{h}}-e^{is\frac{X_j-b}{h}}) \phi_w(s)
\frac{1}{\phi_k(s/h)}
ds,\label{term3*}
\end{eqnarray}
The next lemma establishes asymptotic normality of $W_{nh}(b)$.

\begin{lem}\label{asnor1*}
\begin{equation}
\frac{\sqrt{n}}{h^{(1+\alpha)\lambda+\lambda_0} 
e^{\frac{1}{\mu h^\lambda}}}\,
(W_{nh}(b) - \ex W_{nh}(b))\convd N\Big(0,{A^2\over
2\pi^2}(\mu/\lambda)^{2+2\alpha}\Gamma(\alpha+1)^2\Big).
\end{equation}
\end{lem}

For the term $R_{nh}(a)$ we have the following order bound.
\begin{lem}\label{rllem*}
\begin{equation}\label{rleq*}
R_{nh}(a)-\ex R_{nh}(a)=
o_P\Big({1\over \sqrt{n}}\,h^{(1+\alpha)\lambda+\lambda_0} e^{\frac{1}{\mu h^\lambda}}\Big).
\end{equation}
\end{lem}

Since for $n$ large enough the absolute value of the term
(\ref{term3*}) is smaller than
$$\frac{1}{\pi}(\epsilon/h)^{1-\lambda_0}
(b-a)e^{\frac{1}{\mu}(\frac{\epsilon}{h})^\lambda},$$ a combination of
the two lemmas proves the theorem.

\subsection{Proofs of the lemmas}\label{lemmaproofs}%sec2.3
We shall need the following  analytic lemma in the remaining proofs.

\begin{lem}\label{intexplem}
Assume Condition W.
For $h\to 0$ we have
\begin{equation}\label{intexpansion1}
\int_0^1\phi_w(s)e^{{1\over 2h^2}s^2}ds
\sim Ah^{2+2\alpha}e^{1\over 2h^2}\Gamma(\alpha+1)
\end{equation}
and for $h\to 0$ and $0<\epsilon<1$ and $\beta\geq 0$ fixed we have
\begin{equation}\label{intexpansion2}
\int_\epsilon^1 s^{-\lambda_0}(1-s)^\beta \phi_w(s)
e^{\frac{1}{\mu h^\lambda}s^\lambda}ds
\sim
A ({\mu \over \lambda}h^{\lambda})^{1+\alpha+\beta}
e^{1\over \mu h^\lambda} \Gamma(\alpha+\beta+1)
\end{equation}
\end{lem}

\noindent{\bf Proof}

We only prove (\ref{intexpansion2}). The expansion
(\ref{intexpansion1}) can be shown in the same way.
Notice that 
\begin{eqnarray*}
\lefteqn{
\Big| I_{[0,\frac{1-\epsilon}{h^\lambda}]}(v) (1-h^\lambda v)^{-\lambda_0}
\frac{\phi_w(1-h^\lambda v)}{(h^\lambda v)^\alpha} v^{\alpha+\beta}
e^{\frac{1}{\mu h^\lambda}\{(1-h^\lambda v)^\lambda-1\}   } \Big|}\\
&\leq&\begin{cases}
\epsilon^{-\lambda_0} \Big( 
\sup_{0\leq t\leq 1} \frac{\phi_w(1-t)}{t^\alpha} \Big) v^{\alpha+\beta}
e^{- \frac{1}{\mu}v}\qquad &\text{if $ \lambda \geq 1$},\\
\epsilon^{-\lambda_0} \Big( 
\sup_{0\leq t\leq 1} \frac{\phi_w(1-t)}{t^\alpha} \Big) v^{\alpha+\beta}
e^{- \frac{\lambda}{\mu}v} \qquad &\text{if 
 $0< \lambda < 1$}.
\end{cases}
\end{eqnarray*}

By substituting $s=1-h^\lambda v$ and dominated convergence we get
\begin{eqnarray*}
\lefteqn{
\int_\epsilon^1  s^{-\lambda_0}(1-s)^\beta \phi_w(s) 
e^{ \frac{1}{\mu h^\lambda}s^\lambda } ds}\\
&=&
h^{\lambda(1+\beta)} \int_0^{\frac{1-\epsilon}{h^\lambda}}
(1-h^\lambda v)^{-\lambda_0} \phi_w(1-h^\lambda v)v^\beta
e^{\frac{1}{\mu h^\lambda}(1-h^\lambda v)^\lambda}dv\\
&=&
h^{\lambda(1+\alpha+\beta)}e^{\frac{1}{ \mu h^\lambda}}
\int_0^{\frac{1-\epsilon}{h^\lambda}}(1-h^\lambda v)^{-\lambda_0} 
\frac{\phi_w(1-h^\lambda v)}{(h^\lambda v)^\alpha}v^{\alpha+\beta}
e^{\frac{1}{\mu h^\lambda}\{(1-h^\lambda v)^\lambda-1\}} dv\\
&\sim&
h^{\lambda(1+\alpha+\beta)}e^{\frac{1}{ \mu h^\lambda}}A
\int_0^\infty v^{\alpha +\beta}e^{-\frac{\lambda}{\mu}v} dv\\
&\sim&
(\frac{\mu}{\lambda}h^{\lambda})^{1+\alpha+\beta}e^{\frac{1}{\mu h^\lambda}}
A \Gamma(\alpha+\beta+1)
\end{eqnarray*}

\bigskip

\noindent The next lemma establishes the asymptotic normality.
\begin{lem}\label{asnor}
Let, for a fixed $x$, $U_{nh}(x)$ and $V_{nh}(x)$ be defined by
\begin{eqnarray*}\label{uvdefdef}
U_{nh}(x)
&=&
\frac{1}{\sqrt{n}}\sum_{j=1}^n \Big(\cos\Big(\frac{X_j-x}{h}\Big)
-\ex \cos\Big(\frac{X_j-x}{h}\Big)\Big),\\
V_{nh}(x)
&=&
\frac{1}{\sqrt{n}}\sum_{j=1}^n \Big(\sin\Big(\frac{X_j-x}{h}\Big)
-\ex \sin\Big(\frac{X_j-x}{h}\Big)\Big),
\end{eqnarray*}
then, as $n\to\infty$ and $h\to 0$,
$$
U_{nh}(x) \convd N(0, \frac{1}{2})\quad \mbox{and}\quad 
V_{nh}(x) \convd N(0, \frac{1}{2}).
$$
\end{lem}

\noindent{\bf Proof}

Note that since the density $g$ is the convolution of $f$ and the
density $k$ it is continuous. Write
$$
Y_j=\frac{X_j-x}{h}\bmod 2\pi.
$$
For $0\leq y<2\pi$ we have
\begin{eqnarray*}
\lefteqn{P(Y_j\leq y)=
\sum_{k=-\infty}^\infty P(2k\pi h+x\leq X_j\leq 2k\pi h+yh+x)}\\
&=&
\sum_{k=-\infty}^\infty \int_{2k\pi h+x}^{2k\pi h+yh+x}g(u)du
\sim
\sum_{k=-\infty}^\infty yh\, g(\xi_{k,h}) \\
&=&
\frac{y}{2\pi}\sum_{k=-\infty}^\infty 2\pi h\, g(\xi_{k,h})
\sim
\frac{y}{2\pi}\infint g(u)du=\frac{y}{2\pi}\,,
\end{eqnarray*}
where $\xi_{k,h}$ is a point in the interval $[2k\pi h+x,2k\pi
h+yh+x] \in [2k\pi h+x,2(k+1)\pi h+x]$.  Since $h \to 0$, the last
equivalence follows from a Riemann sum approximation of the integral.

So we have $Y_j\convd U$, where $U$ is uniformly distributed on the
interval $[0, 2\pi]$.  Since the cosine function is bounded and
continuous it then follows that $\ex|\cos Y_j|^p\to \ex |\cos U |^p$,
for all $p>0$. Consequently
$$
\ex  \cos\Big(\frac{X_j-x}{h}\Big) \to \ex \cos U =0
$$
and
$$
\ex  \cos^2\Big(\frac{X_j-x}{h}\Big)  \to \ex \cos^2 U= \frac{1}{2}.
$$

To prove asymptotic normality of $U_{nh}(x)$ note that it is the
normalized sum of an i.i.d. sequence.  It is sufficient to check
whether the Lyapounov condition for the central limit theorem in
Lo\`eve (1977) holds, i.e.  for $\delta>0$,
$$
\frac{\ex| \cos Y_1-\ex  \cos Y_1|^{2+\delta} }{n^{\delta/2} (\var(\cos Y_1))^{1+\delta/2}}
\sim 
\frac{\ex| \cos U |^{2+\delta}} {n^{\delta/2} (\var(\cos U ))^{1+\delta/2}}
\rightarrow 0,
$$
as $n\rightarrow \infty$.
Hence the  condition holds and asymptotic normality of $U_{nh}(x)$ follows.
The proof of asymptotic normality of $V_{nh}(x)$ is similar.

\subsubsection{Proof of Lemma \ref{decomp1}}%lem2.1
Write
\begin{eqnarray*}
\lefteqn{f_{nh}(x)
=\frac{1}{2\pi} \infint e^{-itx} \phi_w(ht) \frac{1}{\phi_k(t)}
\phi_{emp}(t)dt}\\
&=&\frac{1}{2\pi}
 \int_{-\frac{1}{h}}^{\frac{1}{h}} e^{-itx} \phi_w(ht)
e^{\frac{1}{2}t^2}
\phi_{emp}(t)dt\\
&=&
\frac{1}{2\pi h} \int_{-1}^1
e^{-is\frac{x}{h}} \phi_w(s)
e^{\frac{1}{2h^2}s^2}
\phi_{emp}(\frac{s}{h})ds\\
&=&\frac{1}{2\pi nh}
\sum_{j=1}^n
\int_{-1}^1 e^{is\frac{X_j-x}{h}} \phi_w(s)
e^{\frac{1}{2h^2}s^2}ds\\
&=&
\frac{1}{\pi nh}\sum_{j=1}^n
\int_{0}^1 \cos\Big({s\Big(\frac{X_j-x}{h}}\Big)\Big) \phi_w(s)
e^{\frac{1}{2h^2}s^2}ds.
\end{eqnarray*}
Now use some trigonometry to get
\begin{eqnarray}
\lefteqn{\cos\Big(s\Big(\frac{X_j-x}{h}\Big)\Big)=
\cos\Big(\frac{X_j-x}{h}\Big)+\Big(\cos\Big(s\Big(\frac{X_j-x}{h}\Big)\Big)
-\cos\Big(\frac{X_j-x}{h}\Big)\Big)}\nonumber\\
&=&
\cos\Big(\frac{X_j-x}{h}\Big)
-2\sin\Big(\frac{1}{2}(s+1)\Big(\frac{X_j-x}{h}\Big)\Big)
\sin\Big(\frac{1}{2}(s-1)\Big(\frac{X_j-x}{h}\Big)\Big)\label{trig1}\\
&=&
\cos\Big(\frac{X_j-x}{h}\Big) +R_{n,j}(s),\nonumber
\end{eqnarray}
where $R_{n,j}(s)$ is a remainder term satisfying
\begin{equation}\label{Rnbound1}
|R_{n,j}|\leq ( |x| +|X_j|)\Big(\frac{1-s}{h}\Big).
\end{equation}
The bound follows from the inequalities $|\sin x|\leq |x|$.

By Lemma \ref{intexplem} it follows that $f_{nh}(x)$ equals
\begin{eqnarray*}
\lefteqn{ \frac{1}{\pi h}
\int_0^1
\phi_w(s)
e^{\frac{1}{2h^2}s^2}ds
\frac{1}{n}\sum_{j=1}^n \cos\Big(\frac{X_j-x}{h}\Big)
+
\frac{1}{\pi }\frac{1}{n}\sum_{j=1}^n {\tilde R}_{n,j}}\\
&=&
\frac{A}{\pi }
(\Gamma(\alpha+1)+o(1))
h^{1+2\alpha}e^{\frac{1}{2h^2}}\frac{1}{n}
\sum_{j=1}^n \cos\Big(\frac{X_j-x}{h}\Big)
+\frac{1}{n}
\sum_{j=1}^n {\tilde R}_{n,j},
\end{eqnarray*}
where
$$
{\tilde R}_{n,j}=\frac{1}{\pi }\frac{1}{h}\int_\epsilon^1 R_{n,j}(s)
\,\phi_w(s)
e^{\frac{1}{2h^2}s^2}ds.
$$

For the remainder we have, by (\ref{Rnbound1}) and Lemma
\ref{intexplem},
\begin{eqnarray*}
\lefteqn{|\tilde R_{n,j}| \leq \frac{1}{\pi}(|x|+|X_j|)\frac{1}{h}\int_0^1
\Big(\frac{1-s}{h}\Big)
\,\phi_w(s)
e^{\frac{1}{2h^2}s^2}ds}\\
&=&\frac{A}{\pi}(|x|+|X_j|)(\Gamma(\alpha+2)+o(1))h^{2+2\alpha}e^{\frac{1}{
2h^2}}.
\end{eqnarray*}
Hence
$$
\var \tilde R_{n,j} \leq \ex  \tilde R_{n,j}^2 = O\Big(h^{4+4\alpha}e^{\frac{1}{h^2}}\Big)
$$
and
$$
\frac{1}{n}\sum_{j=1}^n ({\tilde R}_{n,j}-\ex {\tilde R}_{n,j})=
O_P\Big(\frac{h^{2+2\alpha}}{ \sqrt{n}}e^{\frac{1}{ 2h^2}}\Big).
$$

Finally we get
\begin{eqnarray*}
\lefteqn{
\frac{\sqrt{n}}{h^{1+2\alpha}e^{\frac{1}{2h^2}}}\,(f_{nh}(x)-\ex f_{nh}(x))}\\
&=&
\frac{A}{\pi}
(\Gamma(\alpha+1)+o(1))U_{nh}(x) +O_P(h),
\end{eqnarray*}
which finishes the proof of Lemma \ref{decomp1}.

\subsubsection{Proof of Lemma \ref{decomp2}}%lem22
Observe  that
\begin{eqnarray}
\lefteqn{F_{nh}(a,b)=
\int_a^b \Big({1\over 2\pi} \int_{-{1\over h}}^{1\over h}
e^{-itx} \phi_w(ht) e^{{1\over 2}t^2} \phi_{emp}(t)\,dt\Big)\,dx}\nonumber\\
&=&
{1\over 2\pi} \int_{-{1\over h}}^{1\over h}
{1 \over {it}} (e^{-ita}-e^{-itb}) \phi_w(ht)
e^{{1\over 2}t^2} \phi_{emp}(t)dt\nonumber\\
&=&
{1\over 2\pi}\int_{-1}^1
{1 \over {is}} (e^{-is{a\over h}}-e^{-is{b\over h}}) \phi_w(s)
e^{{1\over 2h^2}s^2}\phi_{emp}({s\over h})ds\nonumber\\
&=&
{1\over 2\pi n}
\sum_{j=1}^n \int_{-1}^1 {1 \over {is}}
(e^{is{X_j-a\over h}}-e^{is{X_j-b\over h}}) \phi_w(s)
e^{{1\over 2h^2}s^2}ds\nonumber\\
&=&
{1\over \pi n}
\sum_{j=1}^n \int_0^1 {1 \over s}
\Big(\sin\Big(s\Big({X_j-a\over h}\Big)\Big)-\sin\Big(s\Big({X_j-b\over h}\Big)\Big) \phi_w(s)
e^{{1\over 2h^2}s^2}ds\nonumber\\
&=&
{2\over \pi n}
\sum_{j=1}^n \int_0^1 {1 \over s}
\, \cos\Big(s\Big({X_j-(a+b)/2\over h}\Big)\Big)\sin\Big(s\Big({b-a\over 2h}\Big)\Big) \phi_w(s)
e^{{1\over 2h^2}s^2}ds\nonumber\\
&=&
{2\over \pi n}
\sum_{j=1}^n \int_0^\epsilon {1 \over s}
\, \cos\Big(s\Big({X_j-(a+b)/2\over h}\Big)\Big)
\sin\Big(s\Big({b-a\over 2h}\Big)\Big) \phi_w(s)
e^{{1\over 2h^2}s^2}ds\label{int1}\\%2.16
&+&
{2\over \pi n}
\sum_{j=1}^n \int_\epsilon^1 {1 \over s}
\, \cos\Big(s\Big({X_j-(a+b)/2\over h}\Big)\Big)
\sin\Big(s\Big({b-a\over 2h}\Big)\Big)\phi_w(s)
e^{{1\over 2h^2}s^2}ds,\label{int2}%2.17
\end{eqnarray}
where $0<\epsilon<1$.

First note that, since $|\sin x|/|x|\leq 1$, the absolute value of the terms
in the sum (\ref{int1})
are of order
$O((1/h)e^{{1\over 2h^2}\epsilon^2})$.
So the contribution of (\ref{int1}) minus its expectation is of order
\begin{equation}\label{int1bound}
O_P\Big({1\over h\sqrt{n}}\,e^{{1\over 2h^2}\epsilon^2}\Big).
\end{equation}

Next we consider the term (\ref{int2}).
Write $s$ as $1+(s-1)$. Then
\begin{eqnarray}
\lefteqn{\cos\Big(s\Big({X_j-(a+b)/2\over h}\Big)\Big)
\sin\Big(s\Big({b-a\over 2h}\Big)\Big)}\nonumber\\
&=&
\Big\{
\cos\Big({X_j-(a+b)/2\over h}\Big)
\cos\Big((s-1)\Big({X_j-(a+b)/2\over h}\Big)\nonumber\\
&&\qquad-
\sin\Big({X_j-(a+b)/2\over h}\Big)
\sin\Big((s-1)\Big({X_j-(a+b)/2\over h}\Big)
\Big\}\nonumber\\
&\times&
\Big\{
\sin\Big({b-a\over 2h}\Big)\cos\Big((s-1)\Big({b-a\over
2h}\Big)\Big)\\
&&\qquad+
\cos\Big({b-a\over 2h}\Big)\sin\Big((s-1)\Big({b-a\over 2h}\Big)\Big)
\Big\}\label{trig2}\\
&=&
\sin\Big({b-a\over 2h}\Big)\cos\Big({X_j-(a+b)/2\over h}\Big)\nonumber\\
&+&
(s-1)\Big({b-a\over 2h}\Big)\cos\Big({b-a\over 2h}\Big)\cos\Big({X_j-(a+b)/2\over
h}\Big)\nonumber\\
&+&
R_{j,n}(s)\nonumber,
\end{eqnarray}
where the remainder $R_{j,n}$ satisfies
\begin{equation}\label{Rnbound2}
|R_{j,n}(s)|\leq (c_1+c_2|X_j|)\Big\{\Big({{1-s}\over h}\Big)\,
\sin\Big({b-a\over 2h}\Big)
+{{(1-s)^2}\over h^2} +
{{(1-s)^3}\over h^3}\Big\},
\end{equation}
for some positive constants $c_1$ and $c_2$. The bound follows
from the inequalities $|\sin x -x|\leq |x|^2, |1-\cos x|\leq |x|$
and $|1-\cos x|\leq |x|^2$.

By Lemma \ref{intexplem} it follows that the term (\ref{int2}) equals
\begin{eqnarray*}
\lefteqn{{2\over \pi}
\sin\Big({b-a\over 2h}\Big)
\int_\epsilon^1 {1 \over s}
\,\phi_w(s)
e^{{1\over 2h^2}s^2}ds
{1\over n}\sum_{j=1}^n \cos\Big({X_j-(a+b)/2\over h}\Big)}\\
&+&
{2\over \pi}
\Big({b-a\over 2h}\Big)\cos\Big({b-a\over 2h}\Big)
\int_\epsilon^1 {1 \over s}\,(s-1)
\,\phi_w(s)
e^{{1\over 2h^2}s^2}ds
{1\over n}\sum_{j=1}^n \cos\Big({X_j-(a+b)/2\over h}\Big)\\
&+&
{2\over \pi n}\sum_{j=1}^n {\tilde R}_{n,j}\\
&=&
{2A\over \pi}
(\Gamma(\alpha+1)+o(1))\sin\Big({b-a\over 2h}\Big)
h^{2+2\alpha}e^{1\over 2h^2}
{1\over n}\sum_{j=1}^n \cos\Big({X_j-(a+b)/2\over h}\Big)\\
&-&
{A\over \pi}
(b-a)(\Gamma(\alpha+2)+o(1))\cos\Big({b-a\over 2h}\Big)
h^{3+2\alpha}e^{1\over 2h^2}
{1\over n}\sum_{j=1}^n \cos\Big({X_j-(a+b)/2\over h}\Big)\\
&+&
{2\over \pi n}\sum_{j=1}^n {\tilde R}_{n,j},
\end{eqnarray*}
where
$$
{\tilde R}_{n,j}=\int_\epsilon^1 {1 \over s}\,R_{n,j}(s)
\,\phi_w(s)
e^{{1\over 2h^2}s^2}ds.
$$

For the remainder we have, by (\ref{Rnbound2}) and Lemma
\ref{intexplem},
\begin{eqnarray*}
\lefteqn{|\tilde R_{n,j}| \leq}\\
&\leq& (c_1+c_2|X_j|)\int_\epsilon^1 {1 \over s}
\Big(\Big({{1-s}\over h}\Big)\,
\Big|\sin\Big({b-a\over 2h}\Big)\Big|
+{{(1-s)^2}\over h^2} +
{{(1-s)^3}\over h^3}\Big)
\,\phi_w(s)
e^{{1\over 2h^2}s^2}ds\\
&=&
(c_1+c_2|X_j|)A\Big((\Gamma(\alpha+2)+o(1))h^{3+2\alpha}
e^{1\over 2h^2}\Big|\sin\Big({b-a\over 2h}\Big)\Big|\\
&&\qquad +(\Gamma(\alpha+3)+o(1))h^{4+2\alpha}
e^{1\over 2h^2}+(\Gamma(\alpha+4)+o(1))h^{5+2\alpha}
e^{1\over 2h^2}
\Big).\\
\end{eqnarray*}
Hence
$$
\var \tilde R_{n,j} \leq \ex  \tilde R_{n,j}^2 = O\Big(h^{6+4\alpha}e^{1\over
h^2}\Big|\sin\Big({b-a\over 2h}\Big)\Big|^2\Big) +O(h^{8+4\alpha}e^{1\over
h^2} )
$$
and
\begin{equation}
{2\over \pi n}\sum_{j=1}^n ({\tilde R}_{n,j}-\ex {\tilde R}_{n,j})=
O_P\Big({{h^{3+2\alpha}}\over \sqrt{n}}\,e^{1\over
2h^2}\Big|\sin\Big({b-a\over 2h}\Big)\Big| \Big) +
O_P\Big({{h^{4+2\alpha}}\over\sqrt{n}}\,e^{1\over
2h^2}\Big).
\end{equation}

Finally we get
\begin{eqnarray*}
\lefteqn{{\sqrt{n}\over{h^{2+2\alpha}e^{1\over 2h^2}}}\,(F_{nh}(a,b)-\ex F_{nh}(a,b))}\\
&=&
{A\over \pi}
\Big(
2(\Gamma(\alpha+1)+o(1))\sin\Big({b-a\over 2h}\Big)
-(b-a)(\Gamma(\alpha+2)+o(1))\cos\Big({b-a\over 2h}\Big)h
\Big)\\
&&\times
{1\over \sqrt{n}}\sum_{j=1}^n \Big(\cos\Big({X_j-(a+b)/2\over h}\Big)
-\ex \cos\Big({X_j-(a+b)/2\over h}\Big)\Big)\\
&+&
O_P\Big({1\over h^{3+2\alpha}}\,e^{{1\over 2h^2}(\epsilon^2-1)}\Big)
+
O_P\Big(h\Big|\sin\Big({b-a\over 2h}\Big)\Big| \Big) +
O_P\Big(h^2\Big),
\end{eqnarray*}
which finishes the proof of Lemma \ref{decomp2}.

\subsubsection{Proof of Lemma \ref{rllem}}%lem2.4
Note that, by substituting $s=-1+h^2v$, we have
\begin{eqnarray*}
\lefteqn{{1\over 2\pi}
\int_{-1}^{-\epsilon} {1 \over is}
e^{is {X_j-a\over h}  } \phi_w(s)e^{{1\over 2h^2}s^2}ds}\\
&=&
- ih^{2+2\alpha}e^{{1\over 2h^2}}e^{-i {X_j-a\over h} }{1\over 2\pi}
\int_{-\infty}^\infty
\theta_h^+(v)\, e^{ih(X_j-a)v}dv,
\end{eqnarray*}
where
$$
\theta_h^+(v)=I_{[0,{1-\epsilon\over {h^2}}]}(v) {1\over {-1+h^2v}}
{\phi_w(-1+h^2v)\over {(h^2v)^\alpha}}\, v^\alpha e^{-v(1-{1\over 2}h^2v)},
$$
which for $v\geq 0$ the function $\theta_h^+(v)$ converges to
$\theta^+(v)=-Av^\alpha e^{-v}I_{[0,\infty)}(v)$, as $h\to 0$.

Hence
\begin{eqnarray}
\lefteqn{\left|{1\over 2\pi}
\int_{-1}^{-\epsilon} {1 \over is}
e^{is{X_j-a\over h}} \phi_w(s)e^{{1\over 2h^2}s^2}ds\right|}\nonumber\\
&\leq&
\frac{1}{\sqrt{2\pi}}\,h^{2+2\alpha}e^{{1\over 2h^2}}
\left|\frac{1}{\sqrt{2\pi}}\int_{-\infty}^\infty
\theta_h^+(v)\, e^{ih(X_j-a)v}dv\right|\nonumber\\
&=&
\frac{1}{\sqrt{2\pi}}\,h^{2+2\alpha}e^{{1\over 2h^2}}|\Psi_h^+(h(a-X_j))|,
\label{bound1}
\end{eqnarray}
with $\Psi_h^+$ equal to the Fourier transform of $\theta_h^+$.

Similarly we have
\begin{equation}\label{bound2}
\left|{1\over 2\pi}
\int_{\epsilon}^1 {1 \over is}
e^{is{X_j-a\over h}} \phi_w(s)e^{{1\over 2h^2}s^2}ds\right|
\leq
\frac{1}{\sqrt{2\pi}}\,h^{2+2\alpha}e^{{1\over 2h^2}}|\Psi_h^-(h(X_j-a))|,
\end{equation}
where $\Psi_h^-$denotes the Fourier transform of the function $\theta_h^-$ given by
$$
\theta_h^-(v)=I_{[-\frac{1-\epsilon}{h^2},0]}(v)
\frac{1}{1+h^2v}
\frac{\phi_w(1+h^2v)}{(h^2v)^\alpha}\, v^\alpha e^{v(1+{1\over 2}h^2v)},
$$
which converges to $\theta^-(v)=Av^\alpha e^v I_{(-\infty,0]}(v)$, as
$h\to 0$.  

Under the assumption $ah\to-\infty$ we have $h(a-X_j)\to
-\infty$ and $h(X_j-a)\to -\infty$, almost surely, so we can use
elements of the proof of the Riemann Lebesgue lemma, i.e. Theorem
21.39 in Hewitt and Stromberg (1965).  We have
\begin{eqnarray}
\lefteqn{|\Psi_h^+(y)|\leq \frac{1}{2\sqrt{2\pi}}\int_{-\infty}^\infty
|\theta_h^+(v)-\theta_h^+(v-{\pi\over y})|dv}\nonumber\\
&\leq&
\frac{1}{2\sqrt{2\pi}}\int_{-\infty}^\infty
|\theta^+(v)-\theta^+(v-{\pi\over y})|dv
+\frac{1}{\sqrt{2\pi}}\int_{-\infty}^\infty
|\theta_h^+(v)-\theta^+(v)|dv\label{sum}.
\end{eqnarray}
As $|y|\to\infty$ the first term (\ref{sum}) vanishes by the $L_1$
continuity theorem and the second term by dominated convergence with a
majorant similar to the one in the proof of Lemma \ref{intexplem}.
The Fourier transform of the function $\psi_h^-$ can be treated
similarly.

Combining the  bounds (\ref{bound1}) and (\ref{bound2}) we see that
the variance of (\ref{rnh213}) is bounded by
$$
{1\over n}\,h^{4+4\alpha}e^{{1\over h^2}}(\ex |\Psi_h^+(h(a-X_1))|^2+
\ex |\Psi_h^-(h(X_1-a))|^2).
$$
Since $\Psi_h^+$ and $\Psi_h^-$ are bounded functions, $$\ex
|\Psi_h^+(h(a-X_1))|^2+ \ex |\Psi_h^-(h(X_1-a))|^2$$ vanishes by
dominated convergence.  This proves the order bound of the lemma by
the Markov inequality.\hfill$\Box$

\subsubsection{Proof of Lemma \ref{decomp1*}}
Write
\begin{eqnarray}
\lefteqn{f_{nh}(x)
=\frac{1}{2\pi} \infint e^{-itx} \phi_w(ht) \frac{1}{\phi_k(t)}
\phi_{emp}(t)dt}\nonumber\\
&=&
\frac{1}{2\pi}\int_{-\frac{1}{h}}^{\frac{1}{h}} e^{-itx} \phi_w(ht)
\frac{1}{\phi_k(t)}
\phi_{emp}(t)dt\nonumber\\
&=&
\frac{1}{2\pi h} \int_{-1}^1
e^{-is\frac{x}{h}} \phi_w(s)\frac{1}{\phi_k(s/h)}
\phi_{emp}(\frac{s}{h})ds\nonumber\\
&=&\frac{1}{2\pi nh}
\sum_{j=1}^n
\int_{-1}^1 e^{is\frac{X_j-x}{h}} \phi_w(s)
\frac{1}{\phi_k(s/h)}ds\nonumber\\
&=&
\frac{1}{2\pi nh}\sum_{j=1}^n\int_{-\epsilon}^\epsilon
e^{is\frac{X_j-x}{h}}\phi_w(s)\frac{1}{\phi_k(s/h)}ds
\label{epsint}\\
&&\qquad \qquad +
\frac{1}{2\pi nh}\sum_{j=1}^n
\Big(
\int_{-1}^{-\epsilon}+\int_{\epsilon}^1\Big)
e^{is\frac{X_j-x}{h}}\phi_w(s)\frac{1}{\phi_k(s/h)}ds
\label{1int}.
\end{eqnarray}

Note that the integral in (\ref{epsint}) is real valued and that for
$h$ small enough its variance is bounded by
\begin{eqnarray*}
\lefteqn{\var\Big(
\frac{1}{2\pi nh}\sum_{j=1}^n \int_{-\epsilon}^\epsilon
e^{is\frac{X_j-x}{h}}\phi_w(s)\frac{1}{\phi_k(s/h)}ds
\Big)}\\
&\leq&
\frac{1}{4\pi^2 nh^2}\ex
\Big(
\int_{-\epsilon}^\epsilon
e^{is\frac{X_j-x}{h}x}\phi_w(s)\frac{1}{\phi_k(s/h)}ds
\Big)^2\\
&\leq&
\frac{1}{4\pi^2 nh^2 }\Big(
\int_{-\epsilon}^\epsilon \frac{1}{\phi_k(s/h)}ds \Big)^2\\
&\leq&
\frac{1}{4\pi^2 nh^2 } (2\epsilon)^2 \Big(
\sup_{-\epsilon\leq s\leq \epsilon }\frac{1}{|\phi_k(s/h)|}
\Big)^2\\
&\leq&
\frac{2}{\pi^2 }\frac{1}{n} C^{-2} (\epsilon/h)^{2-2\lambda_0}
e^{ \frac{2}{\mu}(\epsilon/h)^{\lambda}}.
\end{eqnarray*}
So the contribution of (\ref{epsint}) minus its expectation is of order
$$
O_P\Big(\frac{1}{\sqrt{n}}(\epsilon/h)^{1-\lambda_0}
e^{ \frac{1}{\mu}(\epsilon/h)^{\lambda}}\Big).
$$
The term (\ref{1int}) can be written as follows.
\begin{eqnarray}
\lefteqn{\frac{1}{2\pi nh}\sum_{j=1}^n 
\Big(\int_{-1}^{-\epsilon}+\int_{\epsilon}^1\Big)
e^{is\frac{X_j-x}{h}}\phi_w(s)\frac{1}{\phi_k(s/h)}ds}\nonumber\\
&=&
\frac{1}{2\pi nhC}\sum_{j=1}^n 
\Big(\int_{-1}^{-\epsilon}+\int_{\epsilon}^1\Big)
e^{is\frac{X_j-x}{h}}\phi_w(s)(\frac{|s|}{h})^{-\lambda_0} 
e^{\frac{1}{\mu}(\frac{|s|}{h})^\lambda}ds\label{versc}\\
&&+
\frac{1}{2\pi nh}\sum_{j=1}^n 
\Big(\int_{-1}^{-\epsilon}+\int_{\epsilon}^1\Big)
e^{is\frac{X_j-x}{h}}\phi_w(s)
\Big(\frac{1}{\phi_k(s/h)}-\frac{1}{C}(\frac{|s|}{h})^{-\lambda_0} 
e^{\frac{1}{\mu}(\frac{|s|}{h})^\lambda}\Big)ds\label{verschil}
\end{eqnarray}
Note that both (\ref{versc}) and (\ref{verschil}) are real.  For (\ref{verschil}) we have 
\begin{eqnarray}
\lefteqn{
\frac{1}{2\pi nhC}\sum_{j=1}^n\Big( \int_{-1}^{-\epsilon}+\int_{\epsilon}^1\Big)
e^{is\frac{X_j-x}{h}}\phi_w(s)(\frac{|s|}{h})^{-\lambda_0} 
e^{\frac{1}{\mu}(\frac{|s|}{h})^\lambda}ds}\nonumber\\
&=&
\frac{1}{\pi nC}h^{\lambda_0-1}\sum_{j=1}^n \int_\epsilon^1
\cos\Big(s\Big({X_j-x\over h}\Big)\Big)\phi_w(s) 
|s|^{-\lambda_0} e^{{1\over \mu h^\lambda}|s|^{\lambda}}ds\label{realint}
\end{eqnarray}

By Lemma \ref{intexplem} and (\ref{trig1}) this equals 
\begin{eqnarray*}
\lefteqn{ 
{1\over \pi}h^{\lambda_0-1}\int_\epsilon^1 \phi_w(s) 
s^{-\lambda_0} e^{{1\over \mu h^\lambda}s^{\lambda}}ds
{1\over n}\sum_{j=1}^n \cos\Big({X_j-x\over h}\Big)+
{1\over n}\sum_{j=1}^n \tilde{R}_{n,j}}\\
&=&
{1\over \pi}A(\Gamma(\alpha+1)+o(1))(\mu / \lambda)^{1+\alpha}
h^{\lambda(1+\alpha)+\lambda_0-1}e^{1\over \mu h^\lambda}
{1\over n}\sum_{j=1}^n \cos\Big({X_j-x\over h}\Big)\\
&+&
{1\over n}\sum_{j=1}^n \tilde{R}_{n,j},
\end{eqnarray*}
where, with $R_{n,j}(s)$ as in (\ref{trig1}), 
$$
\tilde{R}_{n,j}={1\over \pi}h^{\lambda_0-1}\int_\epsilon^1
R_{n,j}(s) \phi_w(s) 
s^{-\lambda_0} e^{{1\over \mu h^\lambda}s^{\lambda}}ds.
$$
For this remainder by (\ref{Rnbound1}) and Lemma \ref{intexplem} we
have
\begin{eqnarray*}
\lefteqn{
|\tilde{R}_{n,j}|\leq {2\over \pi}(|x|+|X_j|)h^{\lambda_0-1}
\int_\epsilon^1
\Big({ 1-s\over h}\Big)\phi_w(s) 
s^{-\lambda_0} e^{{1\over \mu h^\lambda}s^{\lambda}}ds}\\
&=&
{2A\over \pi}(|x|+|X_j|)(\Gamma(\alpha+2)+o(1))
h^{\lambda(2+\alpha)+\lambda_0-2}e^{1\over \mu h^\lambda}.
\end{eqnarray*}
Hence
$$
\var \tilde{R}_{n,j}\leq \ex \tilde{R}_{n,j}^2=O\Big(
h^{2(\lambda(2+\alpha)+\lambda_0-2)}e^{2\over \mu h^\lambda}\Big)
$$
and
$$
\frac{1}{n} \sum_{j=1}^n (\tilde{R}_{n,j}-\ex\tilde{R}_{n,j})=
O_p\Big( {h^{\lambda(2+\alpha)+\lambda_0-2}\over \sqrt{n}}
e^{1\over \mu h^\lambda} \Big).
$$
The variance of (\ref{verschil}) can be bounded by
\begin{eqnarray*}
\lefteqn{ \var\Big(
\frac{1}{2\pi nh}\sum_{j=1}^n 
\Big( \int_{-1}^{-\epsilon} +\int_{\epsilon}^1\Big)
e^{i\frac{X_j-x}{h}s}\phi_w(s)
\Big(\frac{1}{\phi_k(s/h)}-\frac{1}{C}(\frac{|s|}{h})^{-\lambda_0} 
e^{\frac{1}{\mu}(\frac{|s|}{h})^\lambda}\Big)ds\Big)}\\
&\leq&
\frac{1}{4\pi^2 nh^2}\ex
\Big(\Big( \int_{-1}^{-\epsilon} +\int_{\epsilon}^1\Big)
e^{i\frac{X_j-x}{h}s}\phi_w(s)
\Big(\frac{1}{\phi_k(s/h)}-\frac{1}{C}(\frac{|s|}{h})^{-\lambda_0} 
e^{\frac{1}{\mu}(\frac{|s|}{h})^\lambda}\Big)ds \Big)^2\\
&\leq&
\frac{1}{4\pi^2 nh^2C^2}\ex
\Big(\Big( \int_{-1}^{-\epsilon} +\int_{\epsilon}^1\Big)
e^{i\frac{X_j-x}{h}s}\phi_w(s)
(\frac{|s|}{h})^{-\lambda_0} 
e^{\frac{1}{\mu}(\frac{|s|}{h})^\lambda}
\Big| \frac{ C(|s|/h)^{\lambda_0} 
e^{-\frac{1}{\mu}(\frac{|s|}{h})^\lambda}}
{\phi_k(s/h)}-1\Big|ds\Big)^2
\end{eqnarray*}
Since the function 
$$
y \mapsto 
\frac{C|y|^{\lambda_0} e^{\frac{1}{\mu}|y|^\lambda}}{\phi_k(y)}-1$$
is bounded on $\mathbb{R}$ it follows that 
$$
\Big|C \frac{ (|s|/h)^{\lambda_0} 
e^{-\frac{1}{\mu}(\frac{|s|}{h})^\lambda}}
{\phi_k(s/h)}-1\Big|
$$
is bounded and it tends to zero for all fixed $s$ with $|s| \geq
\epsilon$ as $h \to 0$.  So the variance of (\ref{verschil}) is of smaller
order compared to the variance of (\ref{versc}).  This can be shown
by an argument similar to the proofs of lemma \ref{intexplem}.

Finally we get
\begin{eqnarray*}
\lefteqn{
\frac{\sqrt{n}}{h^{\lambda(1+\alpha)+\lambda_0-1}e^{\frac{1}{\mu h^\lambda}}}\,(f_{nh}(x)-\ex f_{nh}(x))}\\
&=&
\frac{A}{\pi}(\frac{\mu}{\lambda})^{1+\alpha}
(\Gamma(\alpha+1)+o(1))U_{nh}(x) +O_P(h^{\lambda-1})
+O_P\Big(\epsilon^{1-\lambda_0}h^{-\lambda(1+\alpha)}
e^{ \frac{1}{\mu h^\lambda}(\epsilon^{\lambda}-1)}\Big),
\end{eqnarray*}
which completes the proof.

\subsubsection{Proof of Lemma \ref{decomp2*}}%lem2.6
Define 
$$
\Delta(X_j,s)=\cos\Big(s\Big({X_j-(a+b)/2\over h}\Big)\Big)
\sin\Big(s\Big({b-a\over 2h}\Big)\Big).
$$
Then, as in the proof of Lemma \ref{decomp2}, $F_{nh}(a,b)$ can be written as follows
\begin{eqnarray}
\lefteqn{F_{nh}(a,b)=
\int_a^b \Big(\frac{1}{2\pi} \int_{-\frac{1}{h}}^{\frac{1}{h}}
e^{-itx} \phi_w(ht)\frac{1}{\phi_k(t)}
\phi_{emp}(t)\,dt\Big)\,dx}\nonumber\\
&=&
{1\over 2\pi n}
\sum_{j=1}^n \int_{-1}^1 {1 \over {is}}
(e^{is{X_j-a\over h}}-e^{is{X_j-b\over h}}) \phi_w(s)
\frac{1}{\phi_k(s/h)}\,ds\nonumber\\
&=&
{1\over 2\pi n}\sum_{j=1}^n \int_{-1}^1 {1 \over s}\,
\Delta(X_j,s)\phi_w(s)\frac{1}{\phi_k(s/h)}\,ds\nonumber\\
&=&
{1\over \pi n}
\sum_{j=1}^n \int_{-\epsilon}^\epsilon {1 \over s}\,
\Delta(X_j,s) \phi_w(s)\frac{1}{\phi_k(s/h)}\,ds\label{int1*}\\
&+&
{1\over \pi n}
\sum_{j=1}^n \Big(\int_{-1}^{-\epsilon}+\int_\epsilon^1\Big) {1 \over s}
\,\Delta(X_j,s) \phi_w(s)\frac{1}{\phi_k(s/h)}\,ds.\label{int2*}
\end{eqnarray}
Note that, since $|\sin x|\leq |x|$, the absolute value of the terms
in the sum (\ref{int1*})
are of order
$O((\epsilon/h)^{1-\lambda_0} e^{\frac{1}{\mu}
(\epsilon/h)^\lambda})$.
So the contribution of (\ref{int1*}) minus its expectation is of order
\begin{equation}\label{int1bound*}
O_P\Big(\frac{1}{\sqrt{n}}(\epsilon/h)^{1-\lambda_0} e^{\frac{1}{\mu}
(\epsilon/h)^\lambda}
\Big).
\end{equation}

Write (\ref{int2*}) as 
\begin{eqnarray}
\lefteqn{ {1\over \pi n}
\sum_{j=1}^n \Big(\int_{-1}^{-\epsilon}+\int_\epsilon^1\Big) {1 \over s}
\,\Delta(X_j,s) \phi_w(s)\frac{1}{\phi_k(s/h)}\,ds}
\nonumber\\
&=&
\frac{1}{\pi nC}
\sum_{j=1}^n \Big(\int_{-1}^{-\epsilon}+\int_\epsilon^1\Big)
{1 \over s}\, \Delta(X_j,s)\phi_w(s) 
(|s|/h)^{-\lambda_0}
e^{\frac{1}{\mu}(\frac{|s|}{h})^\lambda}ds\label{split1}\\
&+&
\frac{1}{\pi n}
\sum_{j=1}^n \Big(\int_{-1}^{-\epsilon}+\int_\epsilon^1\Big)
{1 \over s}\,\Delta(X_j,s)\phi_w(s) 
\Big(\frac{1}{\phi_k(s/h)}-\frac{1}{C}(|s|/h)^{-\lambda_0}
e^{\frac{1}{\mu}(\frac{|s|}{h})^\lambda}\Big)ds.\label{split2}
\end{eqnarray}

Considering the variance of (\ref{split2}) we get
\begin{eqnarray*}
\lefteqn{\var\Big(
\frac{1}{\pi n}
\sum_{j=1}^n \Big(\int_{-1}^{-\epsilon}+\int_\epsilon^1\Big)
{1 \over s}\,\Delta(X_j,s)\phi_w(s) 
\Big(\frac{1}{\phi_k(s/h)}-\frac{1}{C}(|s|/h)^{-\lambda_0}
e^{\frac{1}{\mu}(\frac{|s|}{h})^\lambda}\Big)ds
\Big)}\\
&\leq&\frac{1}{n\pi^2 C^2} \ex\Big(
\Big(\int_{-1}^{-\epsilon}+\int_\epsilon^1\Big)
{1 \over s}\,\Delta(X_j,s)\phi_w(s)(|s|/h)^{-\lambda_0} 
e^{\frac{1}{\mu}(\frac{|s|}{h})^\lambda}
\Big|C\frac{   (|s|/h)^{\lambda_0}e^{-\frac{1}{\mu}(\frac{|s|}{h})^\lambda}}
{\phi_k(s/h)} -1 \Big|ds\Big )^2
\end{eqnarray*}
The expression
$$
\Big|C\frac{   (|s|/h)^{\lambda_0}e^{-\frac{1}{\mu}(\frac{|s|}{h})^\lambda}}
{\phi_k(s/h)} -1 \Big|
$$
is bounded and tends to zero for all $s$ with $|s|\geq \epsilon$ as $h\to 0$.  As in the proof of Lemma \ref{decomp1*} the variance
of the term (\ref{split2}) is of order 
$$
O_P\Big(\frac{1}{\sqrt{n}}(\epsilon/h)^{1-\lambda_0} e^{\frac{1}{\mu}
(\epsilon/h)^\lambda}
\Big).
$$      

Next write $s$ as $1+(s-1)$. Then, as in (\ref{trig2}),
\begin{eqnarray*}
\lefteqn{\Delta(X_j,s)=\cos\Big(s\Big({X_j-(a+b)/2\over h}\Big)\Big)
\sin\Big(s\Big({b-a\over 2h}\Big)\Big)}\\
&=&
\sin\Big({b-a\over 2h}\Big)\cos\Big({X_j-(a+b)/2\over h}\Big)\\
&+&
(s-1)\Big({b-a\over 2h}\Big)\cos\Big({b-a\over 2h}\Big)\cos\Big({X_j-(a+b)/2\over
h}\Big)\\
&+&
R_{j,n}(s),
\end{eqnarray*}
where the remainder $R_{j,n}$ satisfies
\begin{equation}\label{Rnbound2*}
|R_{j,n}(s)|\leq (c_1+c_2|X_j|)\Big\{\Big({{1-s}\over h}\Big)\,
\Big|\sin\Big({b-a\over 2h}\Big)\Big|
+{{(1-s)^2}\over h^2} +
{{(1-s)^3}\over h^3}\Big\},
\end{equation}
for some positive constants $c_1$ and $c_2$.

By Lemma \ref{intexplem} it follows that (\ref{split2}) equals  %lem2.9
\begin{eqnarray*}
\lefteqn{
\frac{1}{\pi}\sin\Big({b-a\over 2h}\Big)
\Big(\int_{-1}^{-\epsilon}+\int_\epsilon^1\Big)
{1 \over s}\,\phi_w(s)(|s|/h)^{-\lambda_0}
e^{\frac{1}{\mu}(\frac{|s|}{h})^\lambda}ds}\\
&&\times
{1\over n}\sum_{j=1}^n \cos\Big({X_j-(a+b)/2\over h}\Big)\\
&+&
{1\over \pi}
\Big({b-a\over 2h}\Big)\cos\Big({b-a\over 2h}\Big)
\Big(\int_{-1}^{-\epsilon}+\int_\epsilon^1\Big)
{1 \over s}\,(s-1)
\,\phi_w(s)(|s|/h)^{-\lambda_0}
e^{\frac{1}{\mu}(\frac{|s|}{h})^\lambda}ds\\
&&\times
{1\over n}\sum_{j=1}^n \cos\Big({X_j-(a+b)/2\over h}\Big)
+
{1\over \pi n}\sum_{j=1}^n {\tilde R}_{n,j}\\
&=&
\frac{2A}{\pi}\Big(\frac{\mu}{\lambda}\Big)^{1+\alpha}
(\Gamma(\alpha+1)+o(1))\sin\Big({b-a\over 2h}\Big)
h^{(1+\alpha)\lambda+\lambda_0}e^{\frac{1}{\mu h^\lambda}}\\
&&\times
{1\over n}\sum_{j=1}^n \cos\Big({X_j-(a+b)/2\over h}\Big)\\
&+&
{A\over \pi}\Big(\frac{\mu}{\lambda}\Big)^{2+\alpha}
(b-a)(\Gamma (\alpha+2)+o(1))\cos\Big({b-a\over 2h}\Big)
h^{(2+\alpha)\lambda+\lambda_0-1}e^{\frac{1}{\mu h^\lambda}}\\
&&\times
{1\over n}\sum_{j=1}^n \cos\Big({X_j-(a+b)/2\over h}\Big)\\
&+&
{1\over \pi n}\sum_{j=1}^n {\tilde R}_{n,j},
\end{eqnarray*}
where
$$
{\tilde R}_{n,j}=\Big(\int_{-1}^{-\epsilon}+\int_\epsilon^1\Big)
 {1 \over s}\,R_{n,j}(s)
\,\phi_w(s)(|s|/h)^{-\lambda_0}
e^{\frac{1}{\mu}(\frac{|s|}{h})^\lambda}ds.
$$

For this remainder we have, by (\ref{Rnbound2*}) and Lemma
\ref{intexplem},
\begin{eqnarray*}
\lefteqn{|\tilde R_{n,j}|
\leq
 (c_1+c_2|X_j|)\Big(\int_{-1}^{-\epsilon}+\int_\epsilon^1\Big)
{1 \over s}
\Big(\Big({{1-s}\over h}\Big)\,
\Big|\sin\Big({a-b\over 2h}\Big)\Big|+}\\
&&+{{(1-s)^2}\over h^2} +
{{(1-s)^3}\over h^3}\Big)
\,\phi_w(s)(|s|/h)^{-\lambda_0}
e^{\frac{1}{\mu}(\frac{|s|}{h})^\lambda}ds\\
&=&
(c_1+c_2|X_j|)A\Big(
(\Gamma(\alpha+2)+o(1))(\frac{\mu}{\lambda})^{2+\alpha}
h^{(2+\alpha)\lambda+\lambda_0-1}e^{\frac{1}{\mu h^\lambda}}
\Big|\sin\Big({a-b\over 2h}\Big)\Big|+\\
&&+
(\Gamma(\alpha+3)+o(1))(\frac{\mu}{\lambda})^{3+\alpha}
h^{(3+\alpha)\lambda+\lambda_0-2}e^{\frac{1}{\mu h^\lambda}}
\Big)+\\
&&\qquad+(\Gamma(\alpha+4)+o(1))(\frac{\mu}{\lambda})^{4+\alpha}
h^{(4+\alpha)\lambda+\lambda_0-3}e^{\frac{1}{\mu h^\lambda}}
\Big).
\end{eqnarray*}
Hence
$$
\var \tilde R_{n,j} \leq \ex  \tilde R_{n,j}^2 = 
O\Big(
h^{2((2+\alpha)\lambda+\lambda_0-1)}e^{\frac{2}{\mu h^\lambda}}
\Big|\sin\Big({a-b\over 2h}\Big)\Big|^2\Big)+O(
h^{2((3+\alpha)\lambda+\lambda_0-2)}e^{\frac{2}{\mu h^\lambda}})
$$
and
\begin{eqnarray}
\lefteqn{
\frac{2}{\pi n}\sum_{j=1}^n ({\tilde R}_{n,j}-\ex {\tilde R}_{n,j})}\nonumber\\
&=&
O_P\Big( \frac{h^{(2+\alpha)\lambda+\lambda_0-1}}{\sqrt{n}}\,
e^{\frac{1}{\mu h^\lambda}}\Big|\sin\Big({a-b\over 2h}\Big)\Big| \Big) +
O_P\Big(\frac{h^{(3+\alpha)\lambda+\lambda_0-2}}{\sqrt{n}}\,
e^{\frac{1}{\mu h^\lambda}}\Big).\label{restop}
\end{eqnarray}

Finally we get
\begin{eqnarray*}
\lefteqn{\frac{\sqrt{n}}{h^{(1+\alpha)\lambda+\lambda_0}
e^{\frac{1}{\mu h^\lambda}}}\,(F_{nh}(a,b)-\ex F_{nh}(a,b))}\\
&=&
\frac{A}{\pi}
\Big(
2\Big(\frac{\mu}{\lambda}\Big)^{1+\alpha}(\Gamma(\alpha+1)+o(1))
\sin\Big({b-a\over 2h}\Big)+\\
&&\qquad+(b-a)\Big(\frac{\mu}{\lambda}\Big)^{2+\alpha}(\Gamma(\alpha+2)+o(1))
\cos\Big({b-a\over 2h}\Big)h^{\lambda-1}
\Big)\\
&&\times
{1\over \sqrt{n}}\sum_{j=1}^n \Big(\cos\Big({X_j-(a+b)/2\over h}\Big)
-\ex \cos\Big({X_j-(a+b)/2\over h}\Big)\Big)\\
&+&
O_P\Big( \epsilon^{1-\lambda_0} h^{-(1+\alpha)^\lambda-1}e^{\frac{1}{\mu h^\lambda}
(\epsilon^\lambda-1)}
\Big)
+
O_P\Big(h^{\lambda-1}\Big|\sin\Big({a-b\over 2h}\Big)\Big| \Big) +
O_P\Big(h^{2\lambda-2}\Big),
\end{eqnarray*}
which finishes the proof of Lemma \ref{decomp2*}.

\subsubsection{Proof of Lemma \ref{rllem*}}
First write
\begin{eqnarray}
\lefteqn{
{1\over 2\pi}\int_{-1}^{-\epsilon} {1 \over is}
e^{is {X_j-a\over h}  } \phi_w(s) \frac{1}{\phi_k(s/h)}ds}\nonumber\\
&=&
{1\over 2\pi}\int_{-1}^{-\epsilon} {1 \over is}
e^{is {X_j-a\over h}  } \phi_w(s)\frac{1}{C}
(|s|/h)^{-\lambda_0}
e^{\frac{1}{\mu}(\frac{|s|}{h})^\lambda}ds\label{sum1}\\
&+&
{1\over 2\pi}\int_{-1}^{-\epsilon} {1 \over is}
e^{is {X_j-a\over h}  } \phi_w(s)\Big(
\frac{1}{\phi_k(s/h)}-\frac{1}{C}(|s|/h)^{-\lambda_0}
e^{\frac{1}{\mu}(\frac{|s|}{h})^\lambda}
\Big)ds.\label{sum2}
\end{eqnarray}

To bound (\ref{sum1}) note that, by substituting $s=-1+h^\lambda v$,
we have
\begin{eqnarray*}
\lefteqn{{1\over 2\pi C}
\int_{-1}^{-\epsilon} {1 \over is}
e^{is {X_j-a\over h}  } \phi_w(s)
|s/h|^{-\lambda_0}e^{\frac{1}{\mu}|s/h|^\lambda}ds}\\
&=&
- ih^{\lambda_0+(1+\alpha)\lambda}e^{\frac{1}{\mu h^\lambda}}
e^{-i {X_j-a\over h} }{1\over 2\pi C}\int_{-\infty}^\infty
\theta_h^+(v)\, e^{ih^{\lambda-1}(X_j-a)v}dv,\\
\end{eqnarray*}
where
$$
\theta_h^+(v)=- I_{[0,{1-\epsilon\over {h^2}}]}(v) (-1+h^\lambda v)^{-\lambda_0-1}
\frac{\phi_w(-1+h^\lambda v)}{(h^\lambda v)^\alpha}\, v^\alpha e^{\frac{1}{\mu h^\lambda}((1-h^\lambda v)^\lambda-1)}.
$$
Note that, for $v\geq 0$ the function $\theta_h^+(v)$ converges to
$\theta^+(v)=-Av^\alpha e^{-\frac{\lambda}{\mu}v}I_{[0,\infty)}(v)$,
as $h\to 0$.

Hence
\begin{eqnarray}
\lefteqn{\left|\frac{1}{2\pi C}
\int_{-1}^{-\epsilon} {1 \over is}
e^{is {X_j-a\over h}  } \phi_w(s)|
s/h|^{-\lambda_0}e^{\frac{1}{\mu}|s/h|^\lambda}ds
\right|}\nonumber\\
&\leq&
\frac{1}{\sqrt{2\pi}C}\,h^{\lambda_0+(1+\alpha)\lambda}e^{\frac{1}{\mu h^\lambda}}
\left|\frac{1}{\sqrt{2\pi}}\int_{-\infty}^\infty
\theta_h^+(v)\, e^{ih^{\lambda-1}(X_j-a)v}dv\right|\nonumber\\
&=&
\frac{1}{\sqrt{2\pi}}\,h^{\lambda_0+(1+\alpha)\lambda}e^{\frac{1}{\mu h^\lambda}}
|\Psi_h^+(h^{\lambda-1}(a-X_j))|,
\label{bound1*}
\end{eqnarray}
with $\Psi_h^+$ the Fourier transform of $\theta_h^+$.  A similar
bound holds for the integral over $[\epsilon,1]$. 
% \label{rllem*}\label{rleq*}
As in the proof of Lemma \ref{rllem} it now follows that the term
(\ref{sum1}) is of the order (\ref{rleq*}).  
 
To deal with (\ref{sum2}) note that the integral can be rewritten as
\begin{eqnarray}
\lefteqn{\frac{1}{\sqrt{2\pi}iC}h^{\lambda_0+(1+\alpha)\lambda} e^{\frac{1}{\mu h^\lambda}}e^{i\frac{X_j-a}{h}} \frac{1}{\sqrt{2\pi}} \int_{-\infty}^{\infty} \eta_h^+(v)e^{ih^{\lambda-1}(X_j-a)v} dv}\nonumber\\
&=& 
\frac{1}{\sqrt{2\pi}iC}h^{\lambda_0+(1+\alpha)\lambda} e^{\frac{1}{\mu h^\lambda}}e^{i\frac{X_j-a}{h}} \tilde{\Psi}_h^+(h^{\lambda-1}(X_j-a)),
\end{eqnarray}
where
$$
\eta_h^+(v)=I_{[0,\frac{1-\epsilon}{h^\lambda}]}(v) (-1+h^\lambda v)^{-1-\lambda_0} \frac{\phi_w(-1+h^\lambda v)}{(h^\lambda v)^\alpha}v^\alpha
e^{\frac{1}{\mu h^\lambda}((1-h^\lambda v)^\lambda-1)}
u\Big(\frac{1-h^\lambda v}{h}\Big),
$$
and $$ u(y)=\frac{C|y|^{\lambda_0}e^{-\frac{1}{\mu}|y|^\lambda}}{\phi_k(y)}-1.
$$
The Fourier transform $\tilde{\Psi}_h^+$ of $\eta_h^+$ can be bounded by
\begin{eqnarray*}
\lefteqn{
|\tilde{\Psi}_h^+(y)|\leq \frac{1}{2\sqrt{2\pi}} \int_{-\infty}^{\infty} |
\eta_h^+(v)-\eta_h^+(v-\frac{\pi}{y})|dv}\\
&\leq&
\frac{1}{\sqrt{2\pi}} \int_{-\infty}^{\infty}|\eta_h^+(v)|\,dv.
\end{eqnarray*}
Since $u(\cdot)$ is bounded and vanishes at plus and minus infinity,
$\eta_h^+(v)$ vanishes as $h\to 0$.  It follows that the term (\ref{sum2}) is also of the order (\ref{rleq*}).  The integral over $[\epsilon,1]$ can be treated similarly.

\subsection{Proof of Theorem \ref{cauchy}}\label{s24}
Note that $e^{-1/h}f_{nh}(x)$ is equal to a sum of independent bounded
random variables that are functions of pairs $(X,Y_h)$ with
$Y_h=(X-x)/h \bmod 2\pi$.  By a similar argument as in the proof of
Lemma \ref{uvdefdef} it follows that $(X,Y_h) \convd (X,U)$ as $h\to
0$, where $X$ and $U$ are independent and $U$ is uniformly distributed
on $[0, 2\pi]$.  By checking the Lyapounov condition we get asymptotic
normality with zero mean and variance equal to
\begin{eqnarray*}
\lefteqn{
\frac{1}{\pi^2}\ex \Big( \frac{1}{1+(X-x)^2}( \cos U +(X-x)\sin U)\Big)^2}\\
&=& \frac{1}{\pi^2}\ex \frac{1}{(1+(X-x)^2)^2 }
( \cos^2 U +(X-x)^2\sin^2 U +2(X-x)\cos U \sin U)\\
&=&
\frac{1}{2\pi^2}\ex \frac{1}{1+(X-x)^2},
\end{eqnarray*}
which completes the proof.

\bigskip

\bigskip

\noindent{\bf \Large Acknowledgment}
The research of the second author has been financed by the Netherlands
Organization for the Advancement of Scientific Research (NWO).

\bigskip

\nocite{*}
\bibliography{as1}

\end{document}